\documentclass[12pt]{article}

\usepackage[english]{babel}
\usepackage{fullpage}
\usepackage{natbib}
\usepackage{fancyhdr}

\normalfont
\usepackage[a4paper,citecolor=blue,linkcolor=blue,colorlinks=true]{hyperref}

\bibpunct{(}{)}{,}{a}{,}{,}

\usepackage{amsmath}
\usepackage{amssymb}
\usepackage{amsfonts}
\usepackage{amsthm}
\usepackage{graphics}
\usepackage{graphicx}
\usepackage[all]{xy}
\usepackage[usenames]{color}
\usepackage[algo2e, ruled]{algorithm2e}
\usepackage{dsfont}
\usepackage{mathrsfs}
\usepackage{fancybox}
\usepackage{float,latexsym,psfrag,epsf,epsfig,amssymb,rotating}
\usepackage{color}
\usepackage{url}
\usepackage{array}
\usepackage{textcomp}

\newcommand{\s}{\mathscr{S}}
\newcommand{\Y}{\mathscr{Y}}
\newcommand{\G}{\mathscr{G}}
\newcommand{\like}{f\left( y \mid \theta\right)}
\newcommand{\pseudo}{f_{\text{pseudo}}\left( y \right)}
\newcommand{\ppseudo}{p_{\text{pseudo}}\left(\theta \mid y \right)}
\newcommand{\CL}{f_{\text{CL}}\left( y \mid\theta \right)}
\newcommand{\post}{p\left( \theta \mid y\right)}
\newcommand{\pCL}{p_{\text{CL}}(\theta\mid y)}

\newcommand{\ivj}{i\stackrel{\G}{\sim} j}
\newcommand{\lvj}{\ell\stackrel{\G}{\sim} j}
\newcommand{\grad}{\mathbf{\nabla}}
\newcommand{\hess}{\mathbf{H}}
\newcommand{\K}{\mathbf{K}}
\newcommand{\E}{\mathbb{E}}
\newcommand{\mapCL}{\widehat{\theta}^{\ast}_{\text{CL}}}
\newcommand{\map}{\widehat{\theta}^{\ast}}
\newcommand{\pCLarg}{p_{\text{CL}}}
\newcommand{\pCLsub}{\overline{p_{\text{CL}}}(\theta\mid y)}
\newcommand{\bfy}{\hbox{\boldmath$y$}}

\begin{document}

\title{\bf Calibration of conditional composite likelihood for Bayesian inference on Gibbs random fields}

\author{
Julien Stoehr\footnote{I3M -- UMR CNRS 5149, 
Universit\'e Montpellier, France.},
\and Nial Friel\footnote{School of Mathematical Sciences and Insight: the national center for data analytics, University College Dublin, Ireland.}
}

\date{\today}
\maketitle

\begin{abstract}
Gibbs random fields play an important role in statistics, however, the resulting likelihood is typically unavailable due to an intractable normalizing constant. Composite likelihoods offer a principled means to construct useful approximations. This paper provides a mean to calibrate the posterior distribution resulting from using a composite likelihood and illustrate its performance in several examples. 
\end{abstract}

\section{Introduction}

Gibbs random fields play an important and varied role in statistics. The autologistic model is used to model the
spatial distribution of binary random variables defined on a lattice or grid \citep{bes74}. The exponential random 
graph model or $p^*$ model is arguably the most popular statistical model in social network analysis \citep{rob:pat:kal:lus07}.
Other application areas include biology, ecology and physics. 

Despite their popularity, Gibbs random fields present considerable difficulties from the point of 
view of parameter estimation, because the likelihood function is typically intractable for all but trivially small graphs. 
One of the earliest approaches to overcome this difficulty is the pseudolikelihood method \citep{bes75}, which replaces 
the joint likelihood function by the product of full-conditional distributions of all nodes. 
It is natural to consider generalizations
which refine pseudolikelihood by considering products of larger collections of variables. The purpose of this paper is to
consider such composite likelihood methods. In particular, we are interested in their use for Bayesian inference. 
\citet{frielproc} focused on a similar problem and studied how the size of the collections of variables influence 
the resulting approximate posterior distribution.
Our main contribution is to present an approach to calibrate the posterior distribution resulting from using a mis-specified
likelihood function.

This paper is organised as follows. Section~\ref{sec:gibbs_ran_fields} outlines a description of Gibbs random fields, and in
particular the autologistic distribution. Composite likelihoods are introduced in Section~\ref{sec:comp_like}. Here we focus
especially on how to formulate conditional composite likelihoods for application to the autologistic model. We also focus on
the issue of calibrating the composite likelihood function for use in a Bayesian context. Section~\ref{sec:examples} illustrates
the performance of the various estimators for simulated data. The paper concludes with some remarks in 
Section~\ref{sec:conclusions}.

\section{Discrete-valued Markov random fields}
\label{sec:gibbs_ran_fields}

A Markov random field $y$ is a family of random variables $y_i$ indexed by a finite set $\s = \{1,\dots,n\}$ of nodes of a graph and taking values from a finite 
state space $\Y$. 
Here the dependence structure is given by an undirected graph $\G$ 
which defines an adjacency relationship between the nodes of $\s$: by definition $i$ and $j$ are adjacent if and only if they are directly connected by an edge 
in the graph $\G$.
The likelihood of $y$ given a vector of parameters $\theta = (\theta_1,\dots,\theta_d)$ is defined as
\begin{equation}
\like \propto \exp(\theta^T s(y)) := q(y|\theta),
\label{eqn:gibbs_like}
\end{equation}
where $s(y) = (s_1(y),\dots,s_d(y))$ is a vector of sufficient statistics. However a major issue arises due to the fact that the normalizing constant in
(\ref{eqn:gibbs_like}), 
\[
 z(\theta) = \sum_{y\in \Y} \exp(\theta^T s(y)),
\]
depends on the parameters $\theta$, and is a summation over all possible realisation of the Gibbs random field. Clearly, $z(\theta)$
is intractable for all but trivially small situations. This poses serious difficulties in terms of estimating the parameter vector
$\theta$. 

One of the earliest approaches to overcome the intractability of (\ref{eqn:gibbs_like}) is the pseudolikelihood method \citep{bes75} 
which approximates the joint distribution of $y$ as the product of full-conditional distributions for each $y_i$,
\[
 \pseudo = \prod_{i=1}^n f(y_i|y_{-i},\theta),
\]
where $y_{-i}$ denotes $y\setminus\{y_i\}$. This approximation has been shown to lead to unreliable estimates of 
$\theta$, see for example, \cite{ryd:tit98}, \cite{fri:pet:rev09}. This is in fact one of the earliest composite likelihood 
approximations, and we will outline work in this area further in Section~\ref{sec:comp_like}. 

The autologistic model, first proposed by \cite{bes72}, is defined on a regular lattice of size $m\times m'$, where
$n=mm'$. It is used to model the spatial distribution of binary variables, taking values $-1,1$. The autologistic model is 
defined in terms of two sufficient statistics,
\[
 s_0(y) = \sum_{i=1}^n y_i, \;\;\; s_1(y) = \sum_{j=1}^n \sum_{\ivj} y_i y_j,
\]
where the notation $\ivj$ means that lattice point $i$ is connected to lattice point $j$ in $\G$. Following this notation, the normalizing constant of 
an autologistic model should be written $z(\theta,\G)$, highlighting that it also depends on a graph of dependency.
Henceforth we assume that the lattice 
points have been indexed from top to bottom in each column and where columns are ordered
from left to right. For example, for a first order neighbourhood model an interior point $y_i$ has neighbours 
$\{y_{i-m}, y_{i-1}, y_{i+1},y_{i+m}\}$. Along the edges of the lattice each point has either $2$ or $3$ neighbours. 
The full-conditional of $y_i$ can be written as
 \begin{equation}
 f(y_i|y_{-i},\theta) \propto \exp ( \theta_0 y_i + \theta_1 y_i( y_{i-m}+y_{i-1}+y_{i+1}+y_{i+m}) ),
\label{eqn:full-con}
\end{equation}
where $y_{-i}$ denotes $y$ excluding $y_i$. As before, the conditional distribution is modified along the edges of the lattice.
The Hammersley-Clifford theorem \citep{bes74} shows the equivalence between the model defined in (\ref{eqn:full-con}) and in 
(\ref{eqn:gibbs_like}). 
The parameter $\theta_0$ controls the relative abundance of $-1$ and $+1$ values and the parameter $\theta_1$ controls the level of spatial aggregation. 
Note that the Ising model is a special case, resulting from $\theta_0=0$.

The auto-models of \cite{bes74} allow variations on the level of dependencies between edges and a potential anisotropy can be introduced on the graph. 
Indeed, consider a set of graphs $\left\lbrace \G_1,\ldots,\G_d\right\rbrace$. Each graph of dependency $\G_k$ induces a summary statistic 
$s_k(y) = \sum_{j=1}^n \sum_{i\stackrel{\G_k}{\sim} j} y_i y_j$. For example, one can consider an anisotropic configuration of a first order 
neighbourhood model: that is edges of $\G_1$ are all the vertical edges of the lattice and edges of $\G_2$ are all the horizontal ones. Then 
an interior point $y_i$ has neighbours $\{y_{i-1}, y_{i+1}\}$ according to $\G_1$ and $\{y_{i-m}, y_{i+m}\}$ according to $\G_2$. Along the edges of 
the lattice each point has either $1$ or $2$ neighbours. This allows to set an interaction strength that differs according to the direction.

\section{Composite likelihoods}
\label{sec:comp_like}

There has been considerable interests in composite likelihoods in the statistics literature. See, \cite{var:reid:firth11} for a 
recent overview. 
Our primary objective is to work with a realisation from an autologistic distribution $y$. According to the previous
section we denote $\s = \{1,\dots,mm'\}$ as an index set for the lattice points. Following \cite{asun10} we
consider a general form of composite likelihood written as
\begin{equation*}
\CL = \prod_{i=1}^C f(y_{A_i}\mid y_{B_i},\theta).
\end{equation*}
Some special cases arise:
\begin{enumerate}
 \item $A_i=A$, $B_i=\emptyset$, $C=1$ corresponds to the full likelihood.
 \item $B_i=\emptyset$ is often termed \textit{marginal composite likelihood}.
 \item $B_i =A\setminus A_i$ is often termed \textit{conditional composite likelihood}.
\end{enumerate}
The focus of this paper is on conditional composite likelihoods, since the autologistic distribution is defined in terms of conditional distributions. 
Note that the pseudolikelihood is a
special case of $3.$ where each $A_i$ is a singleton. We restrict each $A_i$ to be of the same dimension and in particular to correspond to contiguous
square 'blocks' of lattice points of size $k\times k$. In terms of the value of $C$ in case $3.$, an exhaustive 
set of blocks would result in $C=(m-k+1)\times(n-k+1)$. In particular, we allow the collection of blocks $\{A_i\}$
to overlap with one another.

\subsection{Bayesian inference using composite likelihoods}
\label{sec:bayes:cl}

The focus of interest in Bayesian inference is the posterior distribution
\begin{equation}
 p(\theta|y) \propto \like p(\theta).
 \label{eqn:posterior}
\end{equation}
Our proposal here is to replace the true likelihood $\like$ with a conditional composite likelihood,
leading us to focus on the approximated posterior distribution
\[
 \pCL \propto \CL p(\theta). 
\]
Surprisingly, there is very little literature on the use of composite likelihoods in the Bayesian setting, although 
\cite{pauli11} present a discussion on the use of conditional composite likelihoods. Indeed this
paper suggests, following \cite{lin88}, that a composite likelihood should take the general form
\begin{equation}
 \CL = \prod_{i=1}^C f(y_{A_i}\mid y_{B_i},\theta)^{w_i},
\label{eqn:comp_like_approx1}
\end{equation}
where $w_i$ are positive weights. In related work, \citet{frielproc} examined composite likelihood for various block sizes when $w_i = 1$. 
Our paper deals with the issue of calibrating the weights. Before focusing on the tuning of $w_i$, we highlight here the empirical
observation that non-calibrated composite likelihood leads to an approximated posterior distribution with substantially
lower variability than the true posterior distribution, leading to overly precise precision about posterior parameters, see Figure \ref{fig:iso}.

\subsection{Computing full-conditional distributions of $A_i$}

The conditional composite likelihood which we described above relies on evaluating
\begin{equation}
 f(y_{A_i}|y_{-A_i},\theta) = \frac{\exp\left( \theta_0 s_0(y_{A_i}) + s_1(y_{A_i}\mid y_{-A_i}) \right)}{z(\theta,\G,y_{A_i})},
 \label{eqn:gen_recursions}
\end{equation}
where
\[
 s_0(y_{A_i}) = \sum_{j\in{A_i}} y_j, \;\;\; s_1(y_{A_i}|y_{-A_i}) = \sum_{j\in{A_i}}\sum_{\lvj} y_\ell y_j.
\]
Also the normalizing constant now includes the argument $y_{A_i}$ emphasising that it involves a summation over all possible 
realisations of sub-lattices defined on the set $A_i$ and conditioned on the realised $y_{-A_i}$, that is conditioned by all the lattice point of $y_{-A_i}$ 
connected to a lattice point of $y_{A_i}$ by an edge of $\G$. First we describe an approach 
to compute the overall normalizing constant for a lattice, without any conditioning on a boundary. 

Generalised recursions for computing the normalizing constant of general factorisable models such as the autologistic models have
been proposed by \cite{rev:pet04}.
This method applies to autologistic lattices with a small number of rows, up to 
about $20$, and is based on an algebraic simplification due to the reduction in dependence arising from the Markov property. It 
applies to un-normalized likelihoods that can be expressed as a product of factors, each of which is dependent on only a subset 
of the lattice sites. We can write $q(y\mid\theta)$ in factorisable form as
\begin{equation*}
  q(y\mid\theta) = \prod_{i=1}^{n} q_i(\bfy_i\mid\theta),
\end{equation*}
where each factor $q_i$ depends on a subset $\bfy_i = y_i,y_{i+1},\dots,y_{i+m}$ of $y$, where $m$ is 
defined to be the \emph{lag} of the model. We may define each factor as
\begin{equation}
  q_i(\bfy_i,\theta) = \exp\{ \theta_0 y_i + \theta_1 y_i(y_{i+1}+y_{i+m})\}
  \label{eqn:qi}
\end{equation}
for all $i$, except when $i$ corresponds to a lattice point on the last row or last column, in which case $y_{i+1}$ or $y_{i+m}$, 
respectively, drops out of~(\ref{eqn:qi}).

As a result of this factorisation, the summation for the normalizing constant,
\begin{displaymath}
z(\theta,\G) = \sum_{y}\prod_{i=1}^{n} q_i(\bfy_i\mid\theta)
\end{displaymath}
can be represented as
\begin{equation}
z(\theta,\G) = \sum_{y_{n}} q_{n}(\bfy_{n}\mid\theta)  
\dots  \sum_{y_1}
q_1(\bfy_1\mid\theta)
\label{eqn:z_theta}
\end{equation}
which can be computed much more efficiently than the straightforward summation over the $2^n$ possible lattice realisations. Full 
details of a recursive algorithm to compute the above can be found in \cite{rev:pet04}. Note that this algorithm was extended 
in \cite{fri:rue07} to also allow exact draws from $f(y|\theta)$

The minimum lag representation for an autologistic lattice with a first order neighbourhood occurs for $r$ given by the smaller 
of the number of rows or columns in the lattice. Identifying the number of rows with the smaller dimension of the lattice, the 
computation time increases by a factor of two for each additional row, but linearly for additional columns.
It is straightforward to extend this algorithm to allow one to compute the normalizing constant in (\ref{eqn:gen_recursions}),
so that the summation is over the variables $y_{A_i}$ and each factor involves conditioning on the set $y_{-A_i}$.

\section{Bayesian composite likelihood adjustments}

Approximating the true posterior distribution by remplacing the true likelihood by the composite likelihood leads to misspecification in the mean and variance of approximate posterior distribution as shown 
in Figure~\ref{fig:iso}.  The aim of the following Section is to establish identities that links the gradient and the Hessian of the $\log$-posterior for $\theta$ 
to the moments of sufficient statistics with respect to the distribution of the Gibbs random field, whereupon we use these identities to calibrate the weights $w_i$ 
in (\ref{eqn:comp_like_approx1}).

\subsection{An estimation of the gradient and curvature of the posterior distribution}

Using (\ref{eqn:posterior}) as a starting point, we can write the gradient of the $\log$-posterior for $\theta$ as
\[\grad\log\post = s(y) -\grad z(\theta,\G) + \grad\log p(\theta).\]
It is straightforward to show that
\[\grad z(\theta,\G) = \E_{y\mid\theta} s(y),\]

hence the gradient of the $\log$-posterior for $\theta$ can be written as a sum of moments of $s(y)$, namely
\begin{equation}
\grad\log\post = s(y) - \E_{y\mid\theta} s(y) + \grad \log p(\theta).
\label{eqn:gradient-log-post}
\end{equation}
Taking the partial derivatives of the previous expression yields similar identity for the Hessian matrix of the $\log$-posterior for $\theta$,
\begin{equation}
\hess\log\post =  -\K_{y\mid\theta} (s(y)) + \hess\log p(\theta),
\label{eqn:curvature-log-post}
\end{equation}
where $\K_{y\mid\theta} (s(y))$ denotes the covariance matrix of $s(y)$ when $y$ has distribution $\like$.
Similar to (\ref{eqn:gradient-log-post}) and (\ref{eqn:curvature-log-post}), one can express the gradient and Hessian of the log-posterior $\log\pCL$ in terms of moments of the
sufficient statistics.

\subsection{Mean adjustment}

The mean adjustment aims to ensure that the posterior and the approximated posterior distributions 
have the same maximum. Thus, the adjustment here is simply the substitution
\[\pCLsub = \pCLarg(\theta - \theta^{\ast} + \theta^{\ast}_{\text{CL}}\mid y),\]
where $\theta^{\ast}$ and $\theta^{\ast}_{\text{CL}}$, is the maximum \textit{a posteriori} (MAP) of the posterior distribution $\post$ and 
the approximated posterior distribution $\pCL$, respectively.

Addressing the issue of estimation of $\theta^{\ast}$ and $\theta^{\ast}_{\text{CL}}$, we note generally from equation (\ref{eqn:curvature-log-post}) 
that $\log\post$ and $\log\pCL$ are not concave functions. However the Hessian of the log-likelihood is a semi-negative matrix and so is unimodal. 
A reasonable choice of prior, for example with a semi-negative Hessian matrix, will thus lead to a unimodal posterior distribution. 
Care must be taken to ensure convergence of the optimisation algorithms to 
$\theta^{\ast}$ and $\theta^{\ast}_{\text{CL}}$. In particular, we remark that since the approximate posterior distribution is typically very sharp around the 
MAP, as shown in Figure~\ref{fig:iso}, it can be difficult to ensure convergence of gradient based algorithms in reasonable computational time. However, 
in our experiments we have found that using a BFGS algorithm which 
is based on a Hessian matrix approximation using rank-one updates calculated from approximate gradient evaluations,
provided good performance in
our context. Note that in practice, the gradient evaluated in the algorithm is stochastic and based on a 
standard Monte Carlo estimator of the expectation $\mathbb{E}_{y\mid\theta} s_j(y)$.

\begin{algorithm2e}
\caption{MAP estimation }
\label{algo:bfgs}
\KwIn{A lattice $y$}
\KwOut{Estimators $\map$ of $\theta^{\ast}$ and $\mapCL$ of $\theta^{\ast}_{\text{CL}}$}

\medskip

\textbf{estimate} $\theta^{\ast}_{\text{CL}}$ using a BFGS algorithm based on Monte Carlo estimator of $\grad \log \pCL$\;
\textbf{estimate} $\theta^{\ast}$ using a BFGS algorithm based on Monte Carlo estimator of $\grad \log \pCL$ and starting from $\mapCL$\;
\textbf{return} $\map$ and $\mapCL$\;
\end{algorithm2e}

Estimating $\map$ using a random initialization point in BFGS algorithm (see Algorithm \ref{algo:bfgs}) is inefficient. Indeed, estimating 
$\E_{y\mid\theta} s(y)$ is the most cumbersome part of the algorithm and should be done as little as possible. Despite that $\mapCL$ is not equal to $\map$ it is 
usually close and turns out to yield a good initialization to the second BFGS algorithm.

\subsection{Magnitude adjustment}
\label{sub-sec:mag-adj}

The general approach we propose to adjust the covariance of the approximated posterior is to temper the conditional composite 
likelihood with some weights $w_i$ in order to modify its curvature around the mode. We remark that the curvature of a scalar field at its maximum is 
directly linked to the Hessian matrix. Based on that observation, our proposal is to choose $w_i$ such that 
\[\hess\log p(\theta^{\ast}\mid y)=\hess\log\pCLarg (\theta^{\ast}_{\text{CL}}\mid y).\]
Note in our context, there exists no particular reason to weight each blocks differently. Consequently we assume that each block has the same weight and we denote it $w$.

For the sake of simplicity, assume a uniform prior but everything can be easily written for any prior. When $\theta$ is a scalar parameter, writing identity (\ref{eqn:curvature-log-post}) for $\post$ and $\pCL$ yields
\begin{equation}
\label{eqn:weight}
w = \frac{\text{Var}_{y\mid\theta^{\ast}}(s(y))}{\sum_{i=1}^{C} \text{Var}_{y_{A_i}\mid\theta^{\ast}_{\text{CL}}}(s(y_{A_i}\mid y_{-A_i}))}.
\end{equation}
However this approach does not apply when dealing with autologistic models since $\theta\in\mathbb{R}^{d}$ is a vector. We thus have a scalar constraint 
for an equality between the two matrices 
\[\K_{y\mid\theta^{\ast}} (s(y)) = w \sum_{i=1}^C \K_{y\mid\theta^{\ast}_{\text{CL}}} (s(y_{A_i}\mid y_{-A_i})).\]
In Table~\ref{tab-w} we consider some possible identities that are natural to consider in order to choose a reasonable value for $w$. 
The options $w^{(3)}$ and $w^{(4)}$ include only the information contained in the diagonal of each matrix whereas options $w^{(1)}$, $w^{(2)}$ 
and $w^{(5)}$ take advantage of all the information of the covariance matrix. 

\begin{table}[h]
{\caption{\label{tab-w}Weight options for a magnitude adjustment when $\theta\in\mathbb{R}^{d}$}\vspace{2mm}}
\begin{center}
\small{
\begin{tabular}{ll}
\hline\\
$w^{(1)}$: & 
$\left\lbrace
\displaystyle\frac{\det\left[\K_{y\mid\theta^{\ast}} (s(y))\right]}{\det\left[\sum_{i=1}^C \K_{y\mid\theta^{\ast}_{\text{CL}}} (s(y_{A_i}\mid y_{-A_i}))\right]}
\right\rbrace^{1/d}$ \\
&\\
$w^{(2)}$: &
$\displaystyle\frac{1}{d}\text{tr}\left[\K_{y\mid\theta^{\ast}} (s(y))
\left(\sum_{i=1}^C \K_{y\mid\theta^{\ast}_{\text{CL}}} (s(y_{A_i}\mid y_{-A_i}))\right)^{-1}\right]$ \\
&\\
$w^{(3)}$: &
$\displaystyle\frac{1}{d}\cdot\sum_{i=1}^{d}\frac{\text{Var}_{y\mid\theta^{\ast}} (s_i(y))}{\sum_{i=1}^C \text{Var}_{y\mid\theta^{\ast}_{\text{CL}}} (s_j(y_{A_i}\mid y_{-A_i}))}$ \\
&\\
$w^{(4)}$: & 
$\displaystyle\frac{\text{tr}\left[\K_{y\mid\theta^{\ast}} (s(y))\right]}
{\text{tr}\left[\sum_{i=1}^C \K_{y\mid\theta^{\ast}_{\text{CL}}} (s(y_{A_i}\mid y_{-A_i}))\right]}$\\
&\\
$w^{(5)}$: & 
$\displaystyle\sqrt{\frac{\text{tr}\left[\K^2_{y\mid\theta^{\ast}} (s(y))\right]}
{\text{tr}
\left[\left(\sum_{i=1}^C \K_{y\mid\theta^{\ast}_{\text{CL}}} (s(y_{A_i}\mid y_{-A_i}))\right)^2\right]}}$\\
\end{tabular}
}
\end{center}
\end{table}

\subsection{Curvature adjustment}

The adjustment presented in the previous Section only modify the magnitude of the approximated posterior but do not affect its geometry. The weight $w$ similarly 
affects each direction of space parameters and does not take into account a possible modification of the correlation between the variables induced by the use of a 
composite likelihood approximation.
We expect this phenomenon to be particularly important when dealing with models where there is a potential on singletons such as the autologistic model.
Indeed estimation of the abundance parameter and interaction parameter, $\theta_0$ and $\theta_1$, respectively, do not suffer from the same level 
of approximation relating to the independence assumption between blocks. Thus we should move from the general form (\ref{eqn:comp_like_approx1}) 
with a scalar weight on blocks to one involving a matrix of weights.

Following \cite{ribatet2011} in the context of marginal composite likelihood, our strategy is to write 
\[
\like\approx f_{\text{CL}} \left(y\mid
\theta^{\ast}_{\text{CL}}
+W(\theta-\theta^{\ast}_{\text{CL}})\right),
\]
for some constant $d\times d$ matrix $W$. Note the substitution keeps the same maximum but deforms the geometry of the parameter space through the matrix $W$. 

Assume that $W$ is a lower triangular matrix in order to take into account the correlation between the parameter components. The suggestion of \cite{ribatet2011} 
is to choose $W$ in order to satisfy asymptotic properties of maximum composite likelihood estimators when the sample size tends to infinity. Since we 
only have one observation, we do not focus on the asymptotic covariance matrix results 
but rather on the covariance matrix at the estimated MAP. Indeed, we follow the same approach introduced in Section \ref{sub-sec:mag-adj},   
\[\left.\hess\log p(\theta\mid y)\right|_{\theta=\theta^{\ast}}  = \left.\hess\log \pCLarg(
\theta^{\ast}_{\text{CL}}
+W(\theta-\theta^{\ast}_{\text{CL}})\mid y)\right|_{\theta=\theta_{\text{CL}}^{\ast}} ,\]
which is equivalent to 
\[\hess\log p(\theta^{\ast}\mid y) = W^{T}\hess\log \pCLarg(\theta_{\text{CL}}^{\ast}\mid y) W.\]
This leads to a system of equations that can be easily solved. The problem of uniqueness faced by \cite{ribatet2011} due to a Cholesky decomposition does not raise any  issues. Since we have access to a close form of different Hessians through Monte Carlo estimators, solutions are easy to compute and perform the same way.

\section{Examples}
\label{sec:examples}

In this numerical part of the paper, we focus on models defined on a $16\times 16$ lattice and we use exhaustively all $4\times 4$ blocks. For the lattice of this 
dimension the recursions proposed by \cite{fri:rue07} can be used to compute exactly the normalizing constants $z(\theta, \G)$, $z(\theta,\G, y_{A_i})$ and to 
draw exactly from the distribution $\like$ or from the full-conditional distributions of $A_i$ $f(y_{A_i}\mid y_{-A_i}, \theta)$. This exact computation of the 
posterior serves as a ground truth against which to compare with the posterior estimates of $\theta$ using the various composite likelihood estimators. 
Computation was carried out on a desktop PC with six 3.47Ghz processors and with 8Gb of memory. Computing the normalizing constant of each block 
took 0.0004 second of CPU time. One iteration of the BFGS algortihm took 0.09 seconds to estimate the MAP of the composite likelihood and 1 second to estimate 
the MAP of true likelihood. The weight calibration for one dataset took approximately three minutes.
Note that for more realistic situations involving larger lattices, one requires a sampler to draw from the full likelihood such as the 
Swendsen-Wang algorithm \citep{sw87}, however the computational cost of using this algorithm increases dramatically with the size of the lattice. One possible 
alternative is the slice sampler of \citet{mira2001perfect} that provides exact simulations of Ising models.
 
In each experiment, we simulated 100 realisations from the model. For each realisation, we use the BFGS algorithm \ref{algo:bfgs} with an adhoc stopping condition 
to get the estimators $\map$ and $\mapCL$. One iteration of the algorithm is based on a Monte Carlo estimator of either $\E_{y\mid \theta} s(y)$ or 
$\E_{y_{A_i}\mid \theta} s(y_{A_i}\mid y_{-A_{i}}, \theta)$ calculated from 100 exact draws whereas the Monte Carlo estimators of the covariance matrix 
$\K_{y\mid\map} (s(y))$ and $\K_{y\mid\mapCL} (s(y_{A_i}\mid y_{-A_i}))$ are based on $50000$ exact draws. In all experiments we placed uniform priors on $\theta$.

Comparing the posterior $\post$ with the various posterior approximations $\pCL$ requires knowledge of the covariance matrix of $\theta$. We could have used numerical 
integration but we prefered to use a simple MCMC algorithm. In terms of implementation, 7000 iterations were used with a burn in period of 2000 iterations for
each dataset.

\begin{figure}
\centering
\includegraphics[width=.48\textwidth]{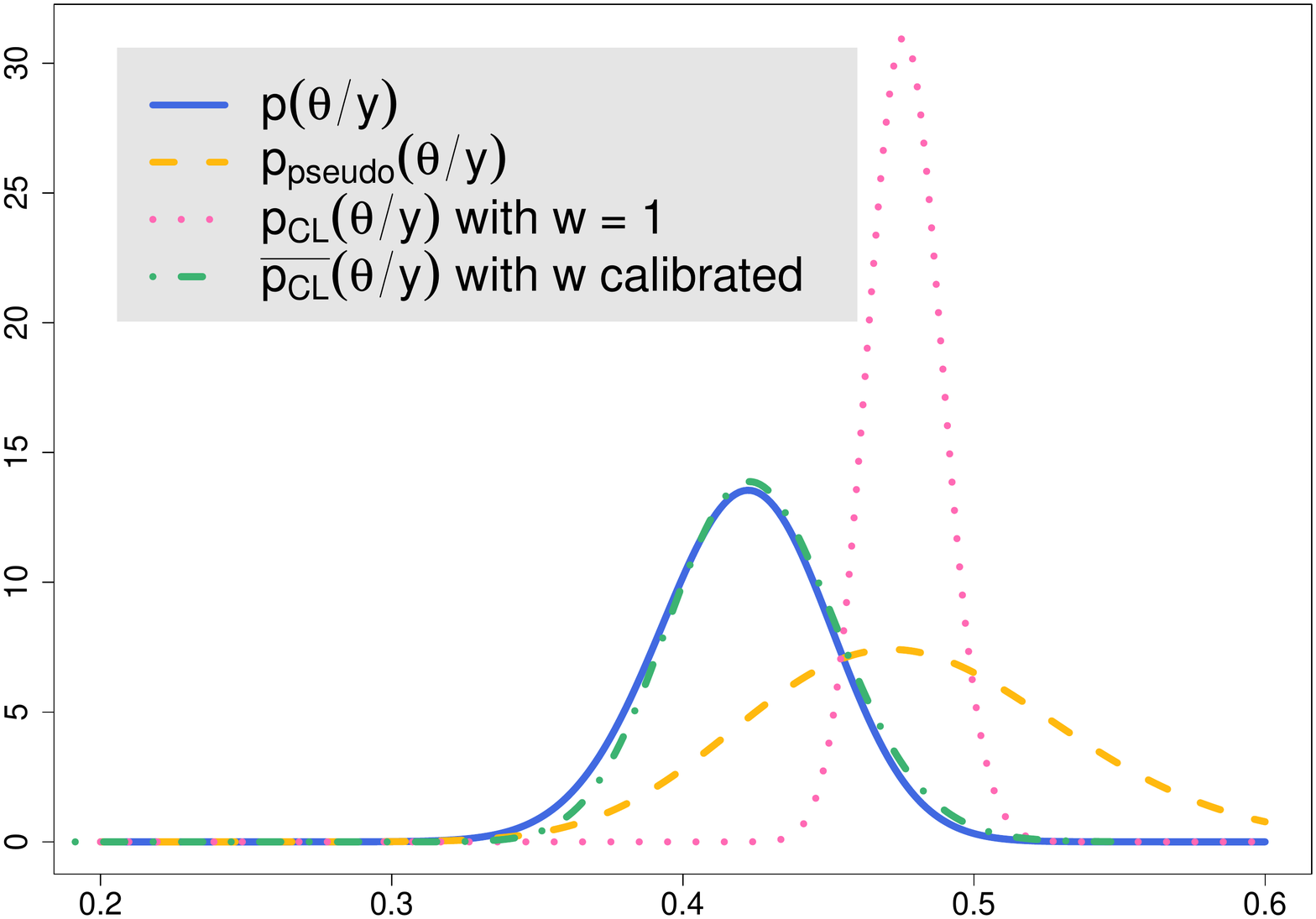} \hfill
\includegraphics[width=.48\textwidth]{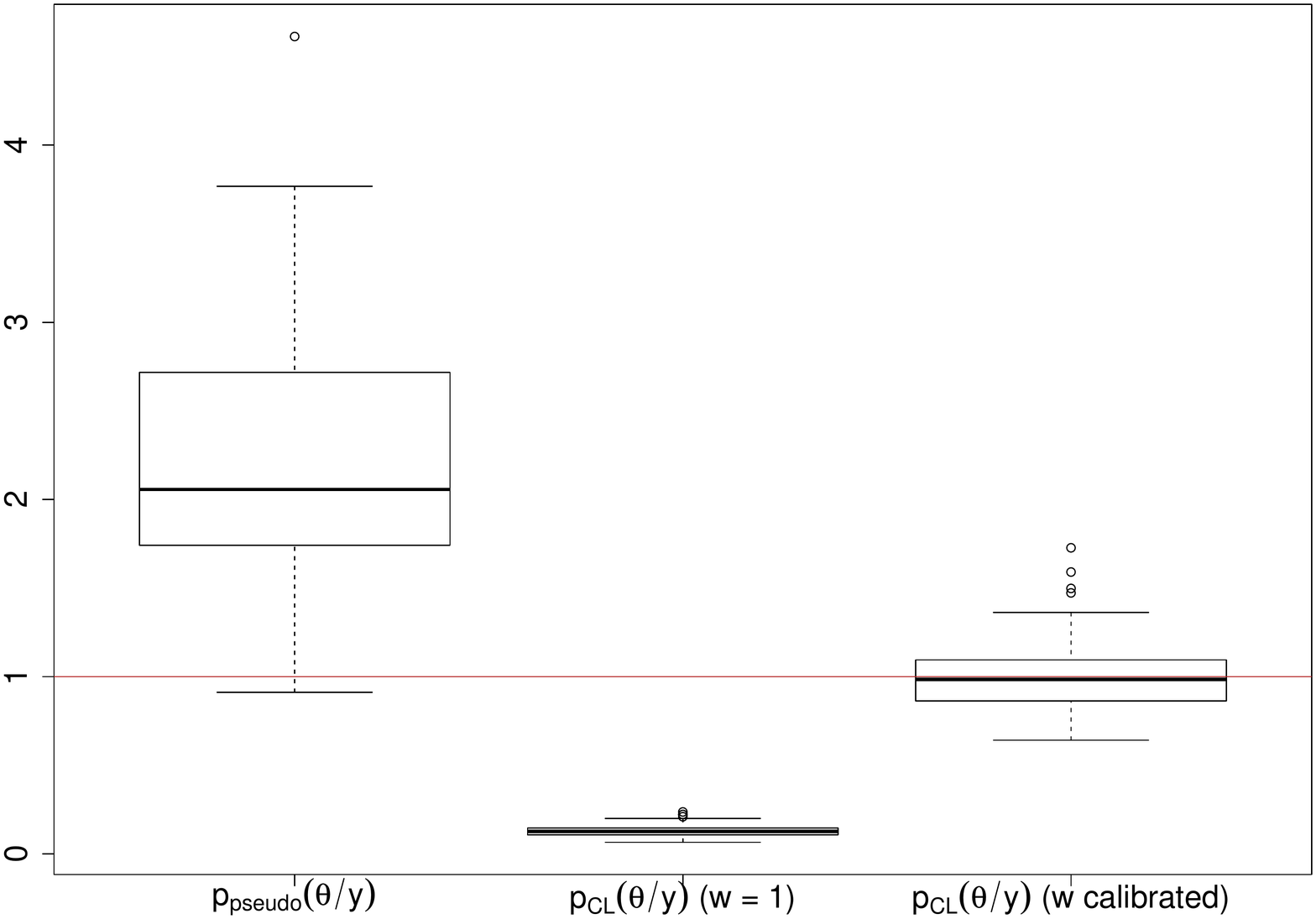}
 \hspace{.2\textwidth} (a) \hspace{.45\textwidth} (b) \hspace{.2\textwidth} ~\\
\caption{\textbf{First experiment results.} (a) Posterior distribution and  posterior distribution approximations for $\theta$ of a first-order Ising model. 
(b) Boxplot displaying the ratio of the variances $\K_{\text{CL}}(\theta)/\K(\theta)$ for 100 realisations of a first-order Ising model.}
\label{fig:iso}
\end{figure}

\bigskip
\textbf{\textit{First experiment}} We considered the special case of a first-order Ising model with a single interaction parameter $\theta = 0.4$, which is close to 
the critical phase transition beyond which all realised lattices takes either value +1 or -1. This parameter setting is the most challenging for the Ising model, since 
realised lattices exhibit strong spatial correlation around this parameter value. Using a fine grid of $\{\theta_k\}$ values, the right hand side of:
\[p(\theta_k\mid y) \propto \frac{q(y\mid\theta_k)}{z(\theta_k)}p(\theta_k), ~k= 1,\ldots,n,\]
can be evaluated exactly. Summing up the right hand side -- using the trapezoidal rule -- yields an estimate of the evidence, $p(y)$, which is the normalizing constant 
for the expression above and which in turn can be used to give a very precise estimate of $p(\theta\mid y)$. The plot so obtained for the posterior and posteriror 
approximations are given by Figure \ref{fig:iso}(a). On this example it should be clear that using an un-calibrated conditional composite likelihood leads to 
considerably underestimated posterior variances. 
But once we perform the mean adjusment and the magnitude adjusment, this provides a very good approximation of the true posterior. In Figure \ref{fig:iso}(b) 
we display the ratio $\K_{\text{CL}}(\theta)/\K(\theta)$, where $\K(\theta)$, respectively 
$\K_{\text{CL}}(\theta)$, denotes the variance of the posterior, respectively the posterior approximation, for $\theta$, based on $100$ realisations of a first-order Ising model.
In view of these results there is no question that the magnitude adjustment (\ref{eqn:weight}) provides an efficient correction of the variance.

Table \ref{tab:iso} confirms this result through evaluation of the relative mean square error, that is $\E\left[(1-\K_{\text{CL}}(\theta)/\K(\theta))^2\right]$,
and the average KL-divergence between the approximated posterior and true posterior distributions based on $100$ 
realisations of a first order.

\begin{table}[h]
\caption{Evaluation of the relative mean square error (RMSE) and the average KL-divergence (AKLD) between the approximated posterior and true posteriror distributions based on $100$ simulations of a first-order Ising model.}
\label{tab:iso}
\begin{center}
  \begin{tabular}{lcc}
  \textbf{COMP. LIKELIHOOD} & \textbf{RMSE} & \textbf{AKLD} \\
    \hline \\
    $\ppseudo$ & $1.96 $ & $0.510 $ 
    \\
    $\pCL$ $(w=1)$ & $0.757 $  & $ 0.337$  
    \\
    $\pCLsub$ ($w$ defined by (\ref{eqn:weight})) & $0.040 $ & $0.010$ 

  \end{tabular}
\end{center}
\end{table}

\begin{figure}

\centering
\includegraphics[width=.48\textwidth]{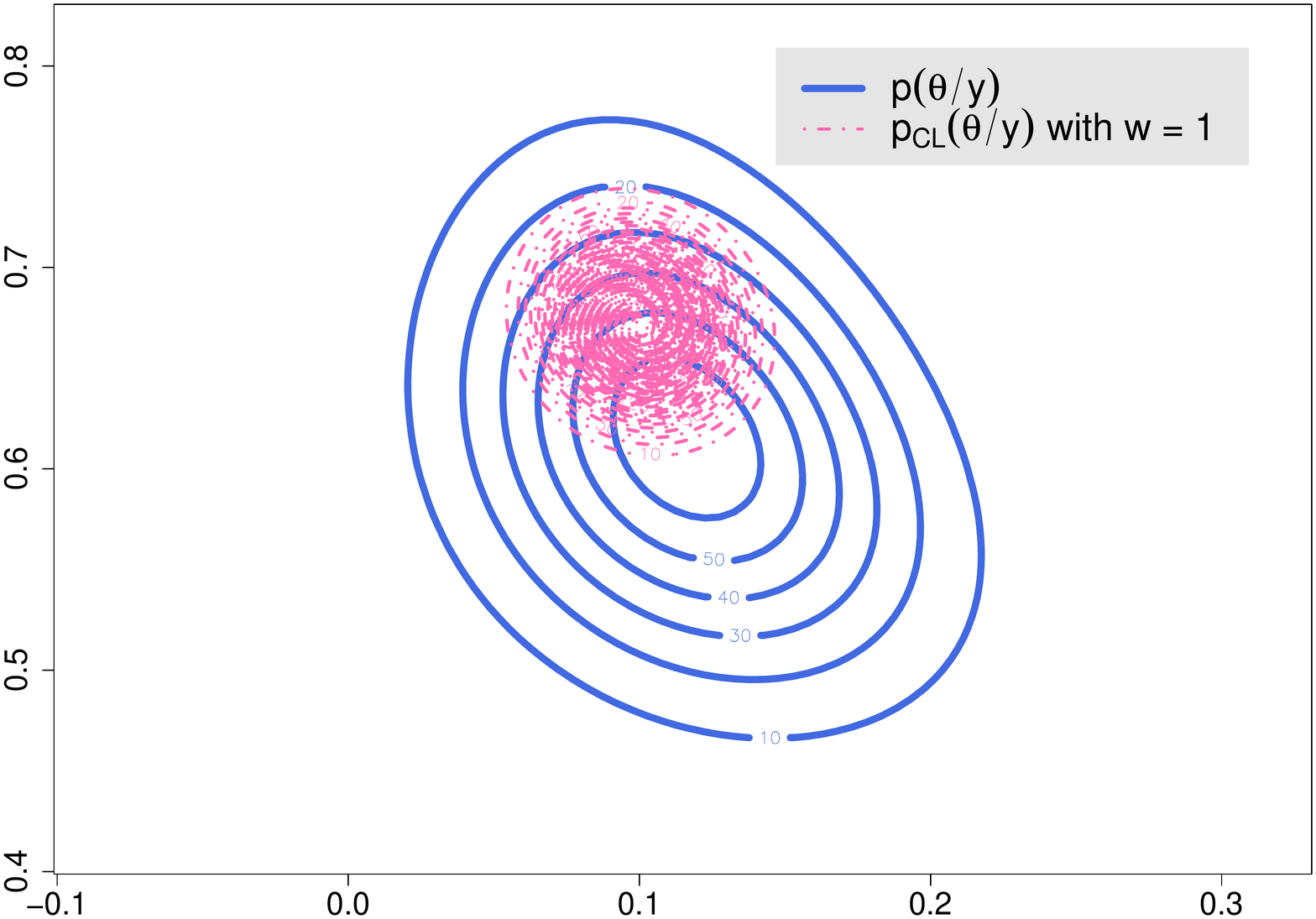}
\hfill
%
\includegraphics[width=.48\textwidth]{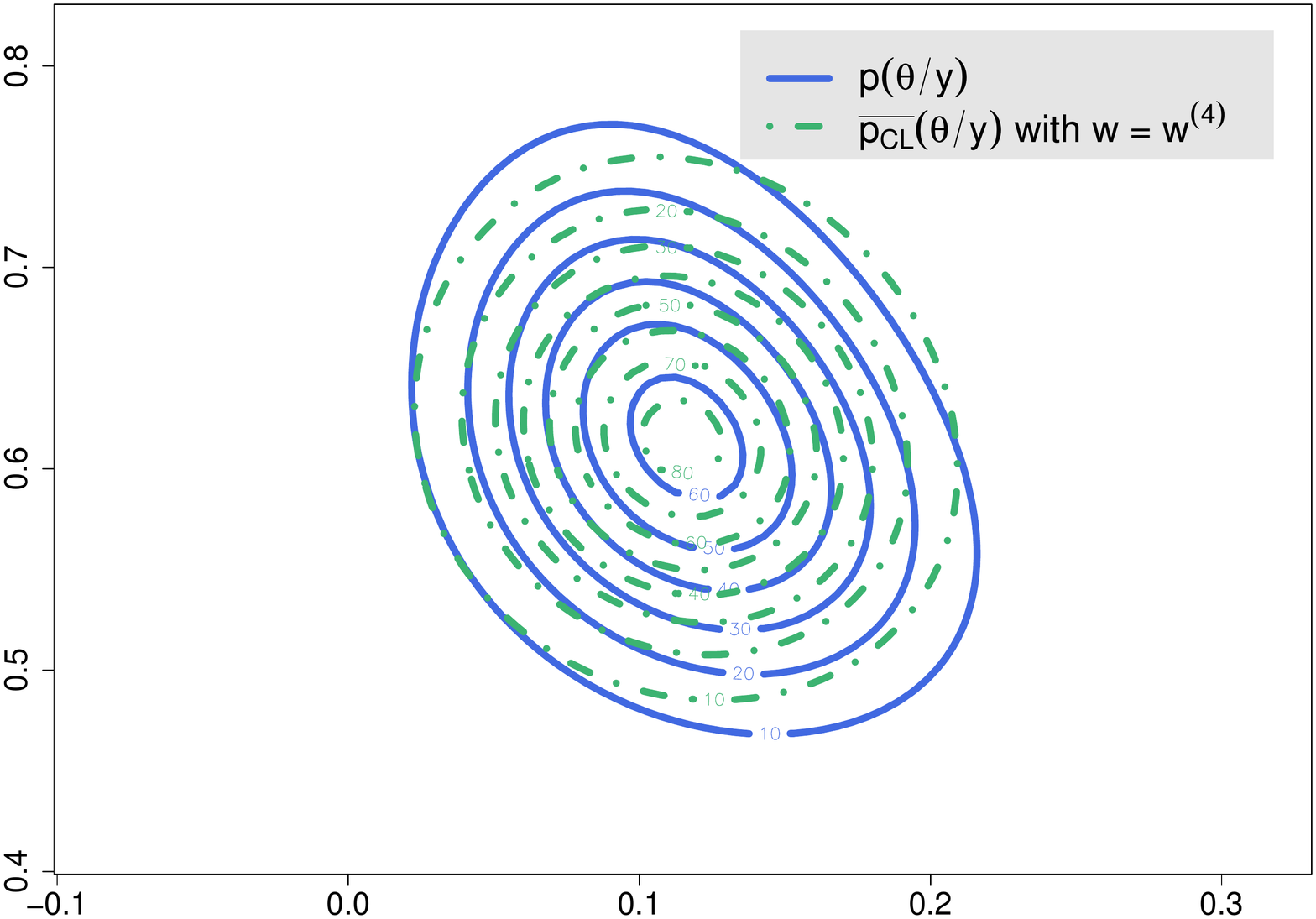}
 \hspace{.2\textwidth} (a) \hspace{.45\textwidth} (b) \hspace{.2\textwidth} ~\\
 \vspace{4mm}
\includegraphics[width=.48\textwidth]{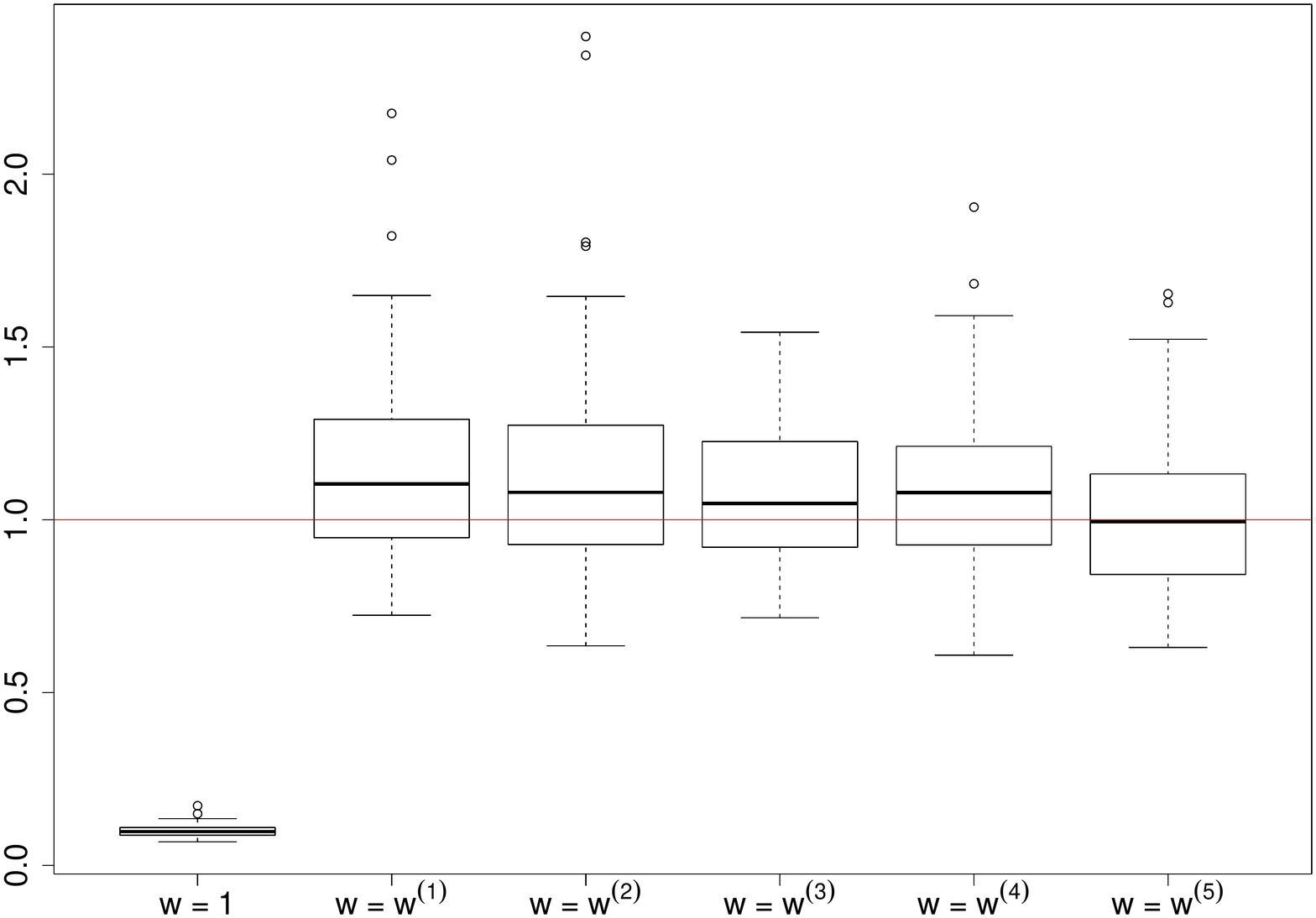}
(c)

\caption{\textbf{Second experiment results.}(a) Posterior distribution and posterior distribution approximation based on the conditional composite likelihood with $w = 1$.
(b) Posterior distribution and posterior distribution approximation based on the conditional composite likelihood with mean and magnitude adjustments ($w = w^{(5)}$).
(c) Boxplot displaying $\Vert\K_{\text{CL}}(\theta)\K^{-1}(\theta)\Vert_{\text{F}}/\sqrt{2}$ for $100$ realisations of an anisotropic first-order Ising model.}
\label{fig:aniso}
\end{figure}

\bigskip
\textbf{\textit{Second experiment}} We were interested in an anisotropic configuration of a first-order Ising model. 
We set $\theta=(0.3,0.5)$. The evidence $p(y)$ is here estimated with an importance sampling method. We drew $1000$ points using a Gaussian law whose moments are 
related to the Monte Carlo estimators of moments of $\theta$. Figure~\ref{fig:aniso}(a) and Figure~\ref{fig:aniso}(b) represent a comparison between the true 
likelihood and the estimates. As for the isotropic case, the mean and the magnitude adjustment allows us to build an accurate approximation of the posterior.
In Figure \ref{fig:aniso}(c) we display boxplots, based on $100$ realisations of an anisotropic first-order Ising model,
of the ratio $\Vert\K_{\text{CL}}(\theta)\K^{-1}(\theta)\Vert_{\text{F}}/\sqrt{2}$, where $\Vert\cdot\Vert_{\text{F}}$ denotes the Frobenius norm. The 
different weight options are almost equivalent in term of variance correction. The weight $w_5$ seems to be the most informative. It should not be a surprise since 
it is based on the Frobenius norm which carries information of the matrix and its singular values.

This conclusion is emphasized in Table \ref{tab:aniso} which presents the relative mean square error 
$\E\left[\Vert1-\K_{\text{CL}}(\theta)\K^{-1}(\theta)\Vert_{\text{F}}^2\right]$ and the average KL-divergence between the approximate and true 
posterior distributions for $100$ realisations of the model.

\begin{table}[h]
  \caption{Evaluation of the relative mean square error (RMSE) the average KL-divergence (AKLD) between the approximated posterior and true posteriror distributions based on $100$ simulations of an anisotropic first-order Ising model.}
  \label{tab:aniso}
\begin{center}
  \begin{tabular}{lcc}
  \textbf{COMP. LIKELIHOOD} & \textbf{RMSE} & \textbf{AKLD}\\
    \hline \\
    $\pCL$ $(w=1)$ & $1.28 $  & $ 2.25$  
    \\
    $\pCLsub$ $(w=w^{(1)})$ & $0.555 $ & $ 0.067$
    \\
    $\pCLsub$ $(w=w^{(2)})$ & $0.583 $ & $ 0.066$ 
    \\
    $\pCLsub$ $(w=w^{(3)})$ & $0.540 $ & $ 0.071$ 
    \\
    $\pCLsub$ $(w=w^{(4)})$ & $0.551 $ & $ 0.061$ 
    \\
    $\pCLsub$ $(w=w^{(5)})$& $0.525 $ & $ 0.079$ 
  \end{tabular}
\end{center}
\end{table}

\begin{figure}[h]
\centering
\includegraphics[width=.48\textwidth]{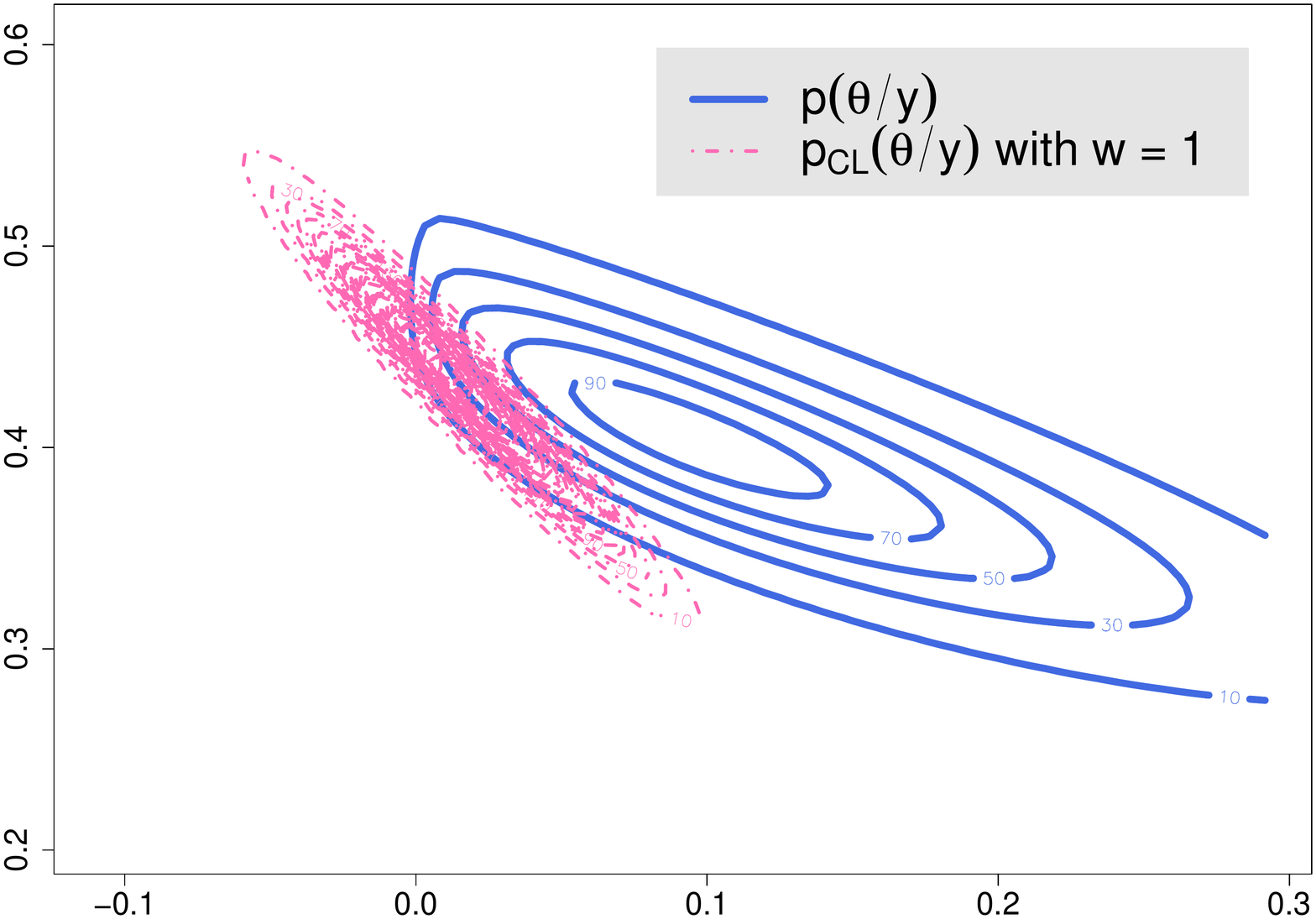}
\hfill
\includegraphics[width=.48\textwidth]{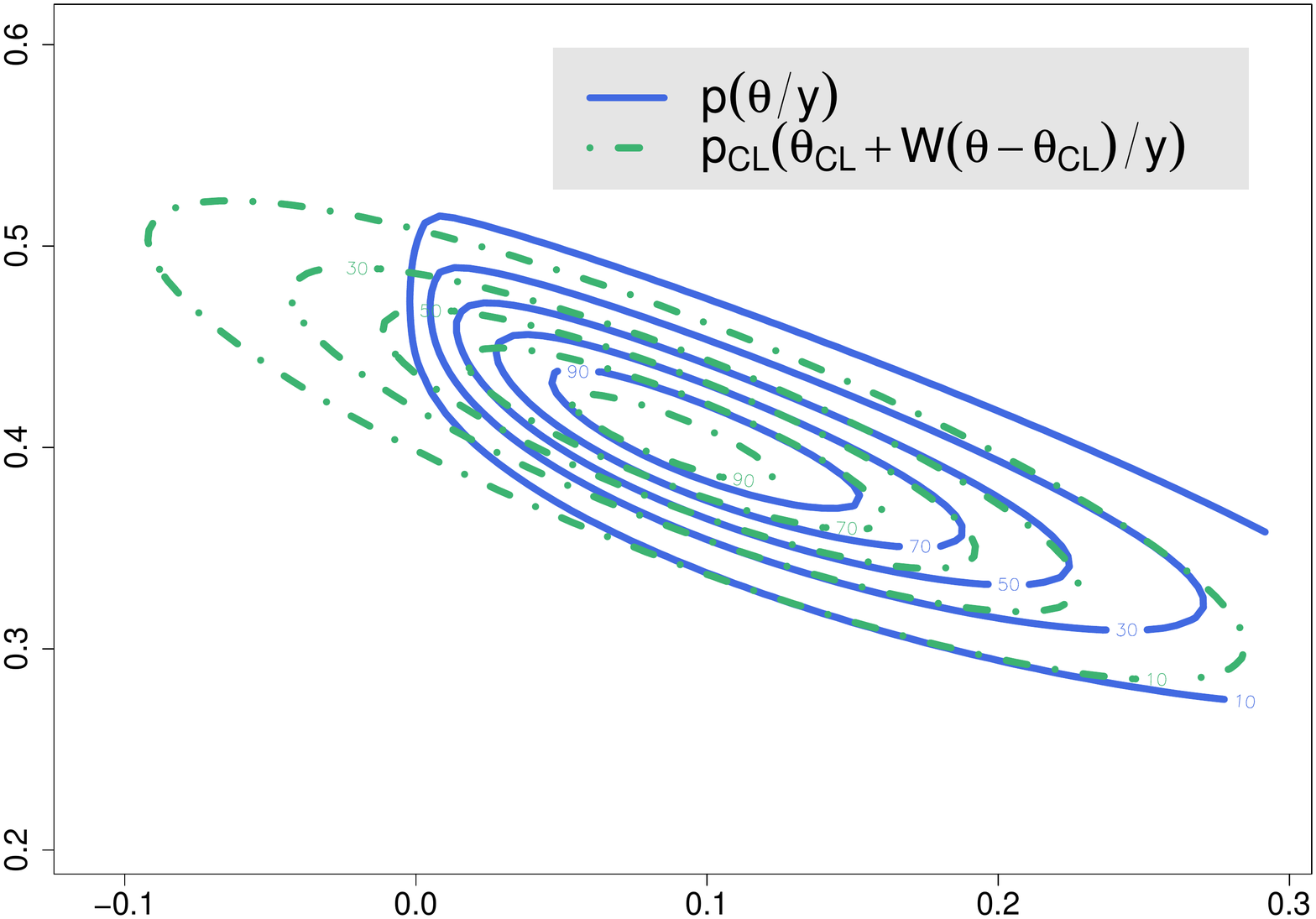}
 \hspace{.2\textwidth} (a) \hspace{.45\textwidth} (b) \hspace{.2\textwidth} ~\\
 \vspace{4mm}
\includegraphics[width=.48\textwidth]{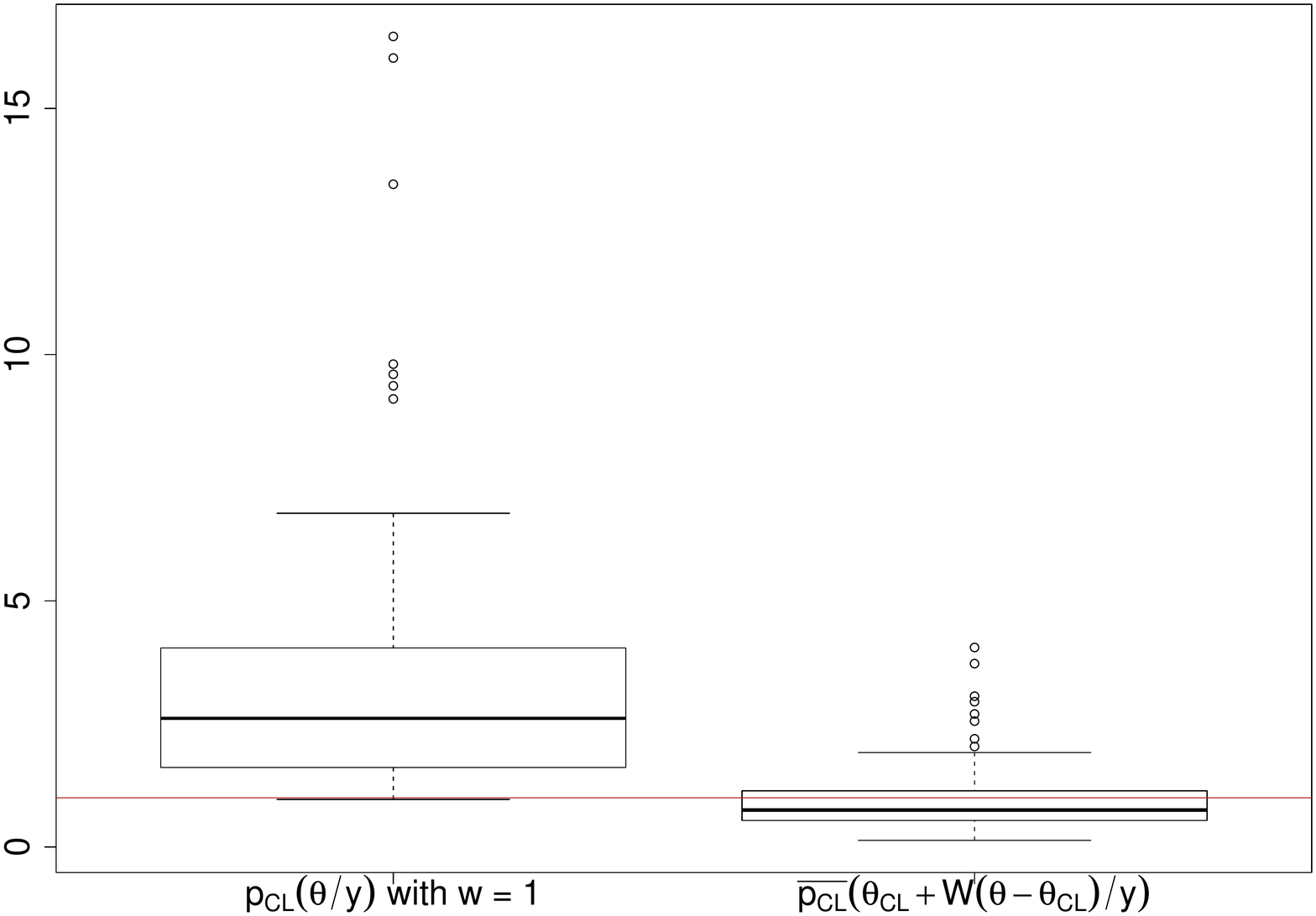}
(c)

\caption{\textbf{Third experiment results.}
(a) Posterior distribution and posterior distribution approximation based on the conditional composite likelihood with $w = 1$.
(b) Posterior distribution and posterior distribution approximation based on the conditional composite likelihood with mean and curvature adjustments.
(c) Boxplot displaying $\Vert\K_{\text{CL}}(\theta)\K^{-1}(\theta)\Vert_{\text{F}}/\sqrt{2}$ for $100$ realisations of a first-order autologistic model.}
\label{fig:abund}
\end{figure}

\bigskip
\textbf{\textit{Third experiment}} Here we focused on an autologistic model with a first-order dependance structure. The abundance parameter was set to $\theta_0 = 0.05$ 
and the interaction parameter to $\theta_1 = 0.4$. The differents implementation settings are exactly the same as for the second experiment. This example 
illustrates how the use of composite likelihood approximation can induce a modification of the geometry of the distribution as shown in Figure \ref{fig:abund}(a). 
Indeed in addition to the mean and variance misspecification the conditional composite likelihood also changes the correlation between the variables. It should be 
evident that a magnitude adjustent would not be fruitful here since it would not affect the correlation. Instead the curvature adjustment manages to do so 
and thus yields a good approximation of the posterior, see Figure~\ref{fig:abund}(b). 
One can object that we do not detect tail of the posterior. But Figure \ref{fig:abund}(c) and Table \ref{tab:abund} show that the adjustment yields an efficient 
correction of the variance.

\begin{table}[h]
   \caption{Evaluation of the relative mean square error (RMSE) and the average KL-divergence (AKLD) between the approximated posterior and true posteriror distributions based on $100$ simulations of a first-order autologistic model.}
  \label{tab:abund}
\begin{center}
  \begin{tabular}{lcc}
 \textbf{COMP. LIKELIHOOD} & \textbf{RMSE} & \textbf{AKLD} \\
    \hline \\
 $\pCL$ $(w=1)$ & $3.44 $  & $ 2.38$  
    \\
    $p_{\text{CL}} \left(
\theta^{\ast}_{\text{CL}}
+W(\theta-\theta^{\ast}_{\text{CL}}) \mid y \right)$ & $0.96 $ & $ 1.89$ 
  \end{tabular}
\end{center}
\end{table}

\section{Conclusion}
\label{sec:conclusions}

This paper has illustrated the important role that conditional composite likelihood approximations can play in the statistical
analysis of Gibbs random fields, and in particular in the Ising and autologistic models in spatial statistics, as
a means to overcoming the intractability of the likelihood function. 
However using composite likelihoods in a Bayesian setting can be problematic, since the resulting approximate posterior distribution
is typically too concentrated and therefore underestimates the posterior mean and variance. 
Our main contribution has been to show how to calibrate the approximate posterior distribution that results from replacing the true 
likelihood with a conditional composite likelihood. Further work will focus on how to extend this framework to Gibbs random fields
with larger number of parameters, such as the exponential random graph model.

\section*{Acknowledgments}
We are grateful to Mathieu Ribatet and the anonymous reviewers of this paper for their helpful comments. 
The Insight Centre for Data Analytics is supported by Science 
Foundation Ireland under Grant Number SFI/12/RC/2289. Nial Friel's research was also supported by an Science Foundation Ireland grant: 12/IP/1424.

\bibliographystyle{abbrvnat}
\bibliography{hmrf}

\end{document}